\documentclass[12pt]{article}
\usepackage[dvips]{graphicx}
\usepackage{amssymb}
\usepackage{amsmath,amsthm}
\usepackage{natbib}
\usepackage{bm}

\newtheorem{theorem}{Theorem}
\newtheorem{lemma}{Lemma}
\newtheorem{corollary}{Corollary}

\newtheorem{proposition}{Proposition}
\newtheorem{assumption}{Assumption}
\newcommand{\R}{{\mathbb{R}}}
\newcommand{\tr}{{\mathrm{tr}}}

\title{Skewness and kurtosis as locally best\\
 invariant tests of normality}
\author{
AKIMICHI TAKEMURA\\
\textit{Graduate School of Information Science and
  Technology}\\
\textit{University of Tokyo}\\
MUNEYA MATSUI\\
\textit{Department of Mathematics}, 
\textit{Keio University}\\
and\\
SATOSHI KURIKI\\
\textit{The Institute of Statistical Mathematics} 
\\
}
\date{August, 2006}
\begin{document}
\maketitle
\begin{abstract}
  Consider testing normality against a one-parameter family of
  univariate distributions containing the normal distribution as the
  boundary, e.g., the family of $t$-distributions or an infinitely
  divisible family with finite variance.  We prove that under mild
  regularity conditions, the sample skewness is the locally best
  invariant (LBI) test of normality against a wide class of asymmetric
  families and the kurtosis is the LBI test against symmetric
  families.  We also discuss non-regular cases such as testing
  normality against the stable family and some related results in the
  multivariate cases.
\end{abstract}


\section{Introduction}
In 1935, E.S.\ Person remarked:
\begin{quote}
  ``\ldots it seems likely that for large samples and 
  when only small departures from normality are in question,  the most
  efficient criteria will be based on the moment  coefficients of the
  sample, e.g.\ on the values of $\sqrt{\beta_1}$ and $\beta_2$.''
\end{quote}
Surprisingly this statement has never been formally proved, although
there exists large literature on testing normality and sampling
distributions of the skewness and the kurtosis.  See \cite{Thod:test:2002} for
a comprehensive survey on tests of normality.  The purpose of this
paper is to give a proof of this statement for fixed sample size ($n
\ge 3$) under general regularity conditions for a wide class of
alternatives, including the normal mixture alternatives and the
infinitely divisible alternatives with finite variance.  Technically all
the necessary ingredients are already given in the literature.
Therefore the merit of this paper is to give a clear statement and a
proof of this basic fact in a unified framework and also to consider
some non-regular cases, in particular testing normality against the
stable family.

In fact ``non-regular'' may not be an appropriate term, because by
considering contamination type alternatives, we see that there are
functional degrees of freedom in constructing an alternative family and
the locally best invariant test against the family.  Therefore by ``small
departure'' we are excluding contamination type departures from normality.
See our discussion at the end of Section \ref{sec:lbi}.

In this paper we are concerned with testing the null hypothesis that the
true distribution belongs to the normal location scale family, against the
alternatives of other location scale families.  We are mainly interested
in invariant testing procedures with respect to the location and the scale
changes of the observations.  In the context of outlier detection,
\cite{Ferg:on:1961} proved that the skewness and the kurtosis are the
locally best invariant tests of normality for slippage type models of
outliers.  In Ferguson's setting, the proportion of outliers can be
substantial but the the amount of slippage tends to zero.  In establishing
the LBI property, \cite{Ferg:on:1961} derived the basic result (see
Proposition \ref{eq:basic-prop} below) on the likelihood ratio of the
maximal invariant under the location-scale transformation.  The same
result was given in Section II.2.2 of
\cite{Hajek:Sidak:1967}.
\cite{Uthoff:1970,Uthoff:1973} used the result to derive the best
invariant tests of normality against some specific alternatives. 
See also Section 3.2 of \cite{Hajek:Sidak:Sen:1999}.
A general result on the likelihood ratio of maximal invariant was given in
\cite{Wijsman:1967,Wijsman:1990} and it led to some important results of
Kariya et al.\ (\cite{Kuwa:Kari:lbi:1991}, \cite{
  Kari:Geor:loca:1994,Kari:Geor:lbi:1995}) in the multivariate setting.

In \cite{Ferg:on:1961}'s setting of outlier detection, if the number of
outliers are distributed according to the binomial distribution, the problem
of outlier detection is logically equivalent to testing normality against
mixture alternatives.  Therefore the LBI property of the skewness and the
kurtosis against mixture alternatives is a straightforward consequence of
\cite{Ferg:on:1961}.  However Ferguson's result has not been interpreted
in this manner.  In this paper we establish the LBI property of the
skewness and the kurtosis in a more general setting and treat the normal
mixture model as an example.

In testing multivariate normality, even if we restrict ourselves to
invariant testing procedures, there is no single LBI test, because 
the maximal invariant moments are multi-dimensional (e.g.\
\cite{Take:maxi:1993}).  Furthermore the invariance can be based on the
full general linear group or the triangular group. This distinction leads
to different results, because the invariance with respect to the triangular
group preserves certain multivariate  one-sided alternatives, whereas 
the general linear group does not.
In Section \ref{sec:multivariate} we discuss these
points in a setting somewhat more general than considered by Kariya and his
coauthors.

The organizations of this paper is as follows. In Section
\ref{sec:lbi} we state our main theorem concerning the locally best
invariant test of normality against one-sided alternatives.  We also
discuss Laplace approximation to the integral in LBI for large sample
sizes $n$.
In Section \ref{sec:mix} we show that our theorem applies in
particular to the normal mixture family and the infinitely divisible
family.  In Section \ref{sec:stable} as an important non-regular case
we consider testing against the stable family.
In Section \ref{sec:profile} we compare locally best invariant test
and tests based on profile likelihood.  Finally in Section
\ref{sec:multivariate} we discuss generalizations of our main theorem to
multivariate cases.

\section{Locally best invariant test of univariate normality}
\label{sec:lbi}

Let 
\[
f_{a,b}(x;\theta)=\frac{1}{b} f\left(\frac{x-a}{b} \; ; \theta\right), \qquad
-\infty < a < \infty, \ b > 0,\  \theta \ge 0 ,
\]
denote a one-parameter family of location-scale densities with the shape
parameter $\theta$.   We simply write
$f(x;\theta)=f_{0,1}(x;\theta)$ for the standard case $(a,b)=(0,1)$.
We assume that $\theta=0$ corresponds to the normal
density 
\[
f(x;0)=\phi(x)=\frac{1}{\sqrt{2\pi}} \exp\Big(-\frac{x^2}{2}\Big).
\]
Based on i.i.d.\ observations $x_1, \ldots, x_n$ from $f_{a,b}(x;\theta)$
we want to test the null hypothesis of normality:
\begin{equation}
\label{eq:one-sided-hypothesis}
H_0:  \theta=0 \quad {\rm vs.}\quad  H_1: \theta>0.
\end{equation}

Here we are testing normality ($\theta=0$) against the one-sided
alternatives
($\theta >0$). 
If we are concerned about heavier tail than the normal as the
alternatives,  this is a natural
setting.  However suppose that we are concerned about asymmetry and we do not know
whether the distribution may be left-skewed or right-skewed under the alternatives.
In this case we should test normality against two-sided alternatives and  then
(\ref{eq:one-sided-hypothesis}) is not a suitable formulation. 
In this paper for simplicity we only consider one-sided alternatives, thus
avoiding the consideration of unbiasedness of tests.

It should be noted that there exists an arbitrariness in choosing a
standard member ($(a,b)=(0,1)$) from a location-scale family.  For the
normal family we usually choose the standard normal density $\phi(x)$ as
the standard member.  Note however that in Section \ref{sec:stable} we
take $N(0,2)$ as the standard member in considering the stable
alternatives for notational convenience.
Given a particular choice of standard members
$f(x;\theta)$, $\theta \ge 0$, we can choose another smooth set of
standard members as
\begin{equation}
\label{eq:another-choice}
f_{a(\theta),b(\theta)}(x;\theta)=
\frac{1}{b(\theta)} f\left(\frac{x-a(\theta)}{b(\theta)} \; ; \theta\right),
\end{equation}
where $a(\theta), b(\theta)$ are smooth function of $\theta$ and
$(a(0),b(0))=(0,1)$.  This arbitrariness does not matter if we use
invariant testing procedures.  However as in the case of normal
mixture distributions in Section \ref{subsec:mixture}, it is sometimes
convenient to resolve this ambiguity in an appropriate manner.
Details on parametrization is discussed in Appendix
\ref{app:parameterization}.

As mentioned above we are primarily interested in invariant testing procedures.  A critical
region $R$ is invariant if
$$
(x_1,\ldots,x_n) \in R  \quad \Leftrightarrow \quad
(a+b x_1,\ldots, a + b x_n) \in R,  \quad
-\infty < \forall a < \infty, \ \forall b >0.
$$
Fix a particular alternative $\theta_1 > 0$.  We state the following basic
result (Theorem b in Section II.2.2 of \cite{Hajek:Sidak:1967}, Section 2
of \cite{Ferg:on:1961}) on the most powerful invariant test against
$\theta_1$.

\begin{proposition}\quad
\label{eq:basic-prop}
The critical region of the most powerful invariant test for testing
$H_0: \theta=0$  against $H_1: \theta = \theta_1 > 0$ is given by
\begin{equation}
\label{eq:lrt}
\frac{\int_0^\infty \int_{-\infty}^\infty \prod_{i=1}^n f(a + b x_i \, ; \theta_1)  b^{n-2} da  db}
{\int_0^\infty \int_{-\infty}^\infty \prod_{i=1}^n f(a + b x_i\, ; 0) b^{n-2} da  db}
> k 
\end{equation}
for some $k > 0$.
\end{proposition}

Note that the values $(x_1,\ldots,x_n)$ can be replaced by any maximal
invariant of the location-scale transformation, since the ratio in
(\ref{eq:lrt}) is invariant.  For our purposes it is most convenient to
replace $x_i$, $i=1,\dots,n$,  by the standardized value 
\begin{equation}
\label{eq:maximal-invariant}
 z_i = \frac{x_i - \bar x}{s}, \quad s^2 =\frac{1}{n} \sum_{j=1}^n (x_j -\bar x)^2.
\end{equation}
Then
$\sum_{i=1}^n z_i=0$ and $\sum_{i=1}^n z_i^2 = n$ and
$$
\prod_{i=1}^n f(a + b
z_i\, ; 0) = 
\frac{1}{(2\pi)^{n/2}}
\exp\Big(-\frac{n(a^2+b^2)}{2}\Big).
$$
Therefore, as in (26) of Section II.2.2 of
\cite{Hajek:Sidak:1967}, the denominator of (\ref{eq:lrt}) becomes the
following constant: 
\begin{eqnarray*}
\int_0^\infty \int_{-\infty}^\infty \prod_{i=1}^n f(a + b
z_i\, ; 0) b^{n-2} da db 
&=& \frac{\Gamma((n-1)/2)}{2 n^{n/2}
   \pi^{(n-1)/2}}.  
\end{eqnarray*}
Since we are considering a fixed sample size $n$, this
constant can be ignored in (\ref{eq:lrt}) and the rejection
region is written as
\begin{equation}
\label{eq:lrt1}
\int_0^\infty \int_{-\infty}^\infty \prod_{i=1}^n f(a + b z_i \, ; \theta_1)
  b^{n-2} da  db > k'.
\end{equation}

We now consider $\theta=\theta_1$ close to $0$.  For a while we proceed formally.
Throughout this paper we assume that 
$l(x;\theta)=\log f(x;\theta)$ is continuously differentiable with respect to
$\theta$ including the boundary $\theta=0$. Then
$$
l(x;\theta)= l(x;0) + l_{\theta}(x;0)  \theta + o(\theta),
$$
where
$$
l_{\theta}(x;\theta)= \frac{\partial}{\partial \theta}\log f(x;\theta)
$$
is the score function.  Therefore
$$
f(x ;\theta)= f(x;0) \exp(l_{\theta}(x;0)\theta + o(\theta))=f(x;0)(1+
l_{\theta}(x;0)\theta)+ o(\theta)
$$
and
\begin{eqnarray*}
\prod_{i=1}^n f(a + b z_i \, ; \theta) 
&=&
\prod_{i=1}^n f(a + b z_i \, ; 0) (1 + \sum_{i=1}^nl_{\theta}(a+bz_i; 0)\theta) +
o(\theta)\\
&=&
\frac{1}{(2\pi)^{n/2}} \exp\Big(-\frac{n(a^2+b^2)}{2}\Big)
(1 + \sum_{i=1}^nl_{\theta}(a+ b z_i; 0)\theta) + o(\theta).
\end{eqnarray*}
It follows that for small $\theta=\theta_1$
 the rejection region (\ref{eq:lrt1}) can be approximately written as
\begin{equation}
\label{eq:lbi}
T(z_1,\ldots,z_n)=\int_0^\infty \int_{-\infty}^\infty 
\sum_{i=1}^nl_{\theta}(a+ b z_i; 0) \exp\Big(-\frac{n(a^2+b^2)}{2}\Big) b^{n-2}
da db > k'{}' .  
\end{equation}

In order to justify the above derivation we assume the following
convenient regularity condition.
\begin{assumption} \quad
\label{ass:1}
For some $\epsilon>0$ 
$$
\int_0^\infty \int_{-\infty}^\infty g(a,b;\epsilon)^n
\exp\Big(-\frac{n(a^2+b^2)}{2}\Big)
b^{n-2} da db < \infty,
$$ 
where
$$
g_n(a,b;\epsilon)=\sup_{|z|\le \sqrt{n}, \;0\le \theta \le \epsilon} 
  \frac{ |\frac{\partial}{\partial \theta} f(a+bz ;
    \theta)|}{f(a+bz ;0 )}.
$$
\end{assumption}  
Under this regularity condition we have the following theorem.

\begin{theorem} \quad  
\label{thm:1}
Under Assumption \ref{ass:1}
the unique rejection region of the locally best invariant
  test of normality $H_0: \theta=0$ vs.\ $H_1: \theta > 0$ 
is given by  (\ref{eq:lbi}), provided that 
$P_0(T(z_1, \ldots, z_n)=k'{}')=0$ under $H_0$.
\end{theorem}

A straightforward proof is given in Appendix \ref{app:proof1}. 
Note that the statement of
this theorem is slightly complicated by the requirement that $P_0(T(z_1,
\ldots, z_n)=k'{}')=0$ under $H_0$.  We need this requirement because if
$P_0(T(z_1,\ldots, z_n)=k'{}')> 0$, in order to maximize the local power
we have to look at $O(\theta^2)$ terms in the expansion of $f(x;\theta)$ 
around $\theta=0$.

A particularly simple result is obtained when $l_{\theta}(x;0)$ is a
polynomial of degree $k$ in $x$.  In this case $l_{\theta}(a+b z_i;0)$
is a polynomial in $a,b$ and $z_i$ and $l_{\theta}(a+b
z_i;0)$ is written as
\begin{equation}
\label{eq:polynomial-score}
l_{\theta}(a+b z_i;0)
= 
p_0(a,b) z_i^k + p_1(a,b)
z_i^{k-1} + \cdots  +   p_k(a,b),
\end{equation}
where  $p_0(a,b),\ldots, p_k(a,b)$ are polynomials in $a$ and $b$.
Denote the standardized $l$-th central moment by
\[
\tilde m_l = \frac{m_l}{s^l}= \frac{1}{n}  \sum_{i=1}^n z_i^l.
\]
Then average of (\ref{eq:polynomial-score}) is written as
\[
\frac{1}{n}
\sum_{i=1}^n l_{\theta}(a+b z_i;0)
= p_0(a,b) \tilde m_k + \dots  
+ p_{k-3}(a,b) \tilde m_3 + p_{k-2}(a,b) +  p_k(a,b).
\]
Furthermore the integral 
$
\int_0^\infty \int_{-\infty}^\infty p_j(a,b)
\exp(-n(a^2+b^2)/2) b^{n-2} da db 
$
can be explicitly evaluated. See Appendix \ref{app:hermite}.
In particular if $l_{\theta}(x;0)$ is a third
degree polynomial, then (\ref{eq:lbi}) is equivalent the standardized
sample skewness of the observations.  Now consider the case that $l_{\theta}(x;0)$
is a fourth degree polynomial without odd degree terms.  Then 
\[
\int_{-\infty}^\infty a^{2l+1} \exp\Big(- \frac{na^2}{2}\Big)da =0 
\]
implies that
$\int_{-\infty}^\infty p_{k-3}(a,b)da=0$ in (\ref{eq:polynomial-score}).
Therefore (\ref{eq:lbi}) is equivalent the standardized sample kurtosis.
We now have the following corollary.

\begin{corollary}\quad
\label{co:1}
  Assume the same regularity condition as in Theorem \ref{thm:1}.
If the score function $l_{\theta}(x;0)$ is a third degree
  polynomial in $x$, then the locally best invariant test of
  normality is given by the standardized sample skewness. 
If  $l_{\theta}(x;0)$ is a fourth degree
  polynomial in $x$  without odd degree terms, then
the locally best invariant test of   normality is given by the
standardized  sample kurtosis.
\end{corollary}

In the next section we show that in two important cases,
$l_{\theta}(x;0)$ is a third degree polynomial for asymmetric
alternatives and is a fourth degree polynomial in $x$ without odd
degree terms for symmetric alternatives.

For general score function the integral (\ref{eq:lbi}) may not be easy
to evaluate.  Although in this paper we are considering fixed $n$, we
here discuss Laplace approximation to the integral (\ref{eq:lbi}) for
large $n$.  Let $A$ denote a random variable having the distribution
$N(0,1/n)$ and let $B$ denote the random variable such that
$B/\sqrt{n}$ has the $\chi$-distribution with $n-1$ degrees of freedom.
Then as $n\rightarrow\infty$, $(A,B)$ converges to $(0,1)$ in
distribution (or equivalently in probability). Note that except for the
normalizing constant, the integral in (\ref{eq:lbi}) can be written as
$E[\sum_{i=1}^n l_\theta(A + B z_i;0)]$.
Under  mild regularity conditions, for large $n$, this expectation is simply
approximated by putting $(A,B)=(0,1)$: 
\begin{equation}
\label{eq:app-LBI}
\tilde T(z_1,\dots,z_n) = \sum_{i=1}^n l_\theta(z_i;0)
\end{equation}
It is easily shown that this is in fact the Laplace approximation
(e.g.\ \cite{Bleistein:Handelsman}) to the integral in (\ref{eq:lbi}).
We call $\tilde T$ approximate LBI for testing normality.  Under mild
regularity conditions, the approximate LBI and the LBI should be
asymptotically equivalent.

In Appendix A.1 of \cite{Kuriki:Takemura:2001} it is shown that the
test based on the $k$-th standardized sample cumulant is
asymptotically equivalent to the test based on $\sum_{i=1}^n
H_k(z_i)$, where $H_k$ is the $k$-th Hermite polynomial.  We see that
the $k$-th standardized sample cumulant is characterized as an
approximate LBI for the case that the score function is given by
$H_k(x)$.  See a further discussion in Section \ref{sec:profile}.
When $n$ is not too large, we may consider evaluating $E[\sum_{i=1}^n
l_\theta(A + B z_i;0)]$ by numerical integration or by Monte Carlo
sampling.

For the rest of this section we make several remarks on the above
results.  In the location-scale transformation $x_i \mapsto a+ b x_i$
we might allow $b\neq 0$ to be negative.  
The maximal invariant is $\bm{z}=(z_1,\ldots,z_n)'$ with $\bm{z}$
identified with $-\bm{z}$, or more compactly it is $\bm{z}\bm{z}'$.
Then an invariant critical region can not depend on a sign preserving
function $\psi$ of $\bm{z}$ (i.e.\ $\psi(-\bm{z})=-\psi(\bm{z}))$.  In
particular it can not depend on the skewness $m_3$ itself, although it
can depend on $|m_3|$.  In the univariate case, allowing $b<0$ is
somewhat unnatural and we have so far only considered $b>0$.  However
in the multivariate case the invariance with respect to the full
general linear group corresponds to allowing $b<0$ in the univariate
case.  We discuss this point further in Section
\ref{sec:multivariate}.

Let $g(x)$ be a probability density.  By an $\epsilon$-contamination
alternative we mean a density of the form $$
f(x;\epsilon)=
(1-\epsilon)\phi(x) + \epsilon g(x)= \phi(x) + \epsilon
(g(x)-\phi(x)).  $$
Letting $\theta = \epsilon$, we see $$
l_{\theta}(x;0)= \frac{g(x)}{\phi(x)} - 1.  $$
Therefore as long as
$g(x)=\phi(x)(1+l_{\theta}(x;0))$ is a probability density, we can
construct a one-parameter contamination family of alternatives such
that $T(z_1,\ldots,z_n)$ in (\ref{eq:lbi}) is the LBI with this score
function $l_{\theta}(x;0)$. By ``small departures from normality''
\cite{Pearson:1935} probably did not have a contamination alternative
in mind.  In our setting the sample size $n$ is fixed.  If $\epsilon$
is much smaller than $1/n$, we actually have no observation from
$g(x)$ with probability close to 1.  In this sense a contamination
family seems to possess certain non-regularity as a family containing
the normal distribution.

\section{Normal mixture family and infinitely divisible family of distributions}
\label{sec:mix}

In this section we discuss two general classes of alternatives such that
the score function at the normal distribution is a polynomial and Corollary
\ref{co:1} is applicable.  The first is the normal mixture family and the
second is the infinitely divisible family with finite variance.

\subsection{Normal mixture family}
\label{subsec:mixture}

Suppose that the mean $\mu$ and the variance $\sigma^2$ of the normal
distribution $N(\mu,\sigma^2)$ has the prior
distribution $g(\mu,\sigma^2;\theta)$, \ $\theta\ge 0$, such that
$g$ degenerates to the point mass at $(0,1)$ as $\theta \rightarrow 0$.
For simplicity write $\tau=1/\sigma^2-1$.  Then as $\theta\rightarrow 0$,
both $\mu$ and $\tau$ converge to $0$ in distribution.
The marginal density is given by
$$
f(x;\theta) = \int_{-1}^\infty \int_{-\infty}^\infty \frac{1}{\sqrt{2\pi}}
\exp(- (\tau+1) \frac{(x-\mu)^2}{2}) h(\mu,\tau;\theta) d\mu d\tau, 
$$
where $h(\mu,\tau;\theta)=(\tau+1)^2 g(\mu,1/(1+\tau);\theta)$.  
Consider the expansion
\begin{eqnarray}
\exp(- (\tau+1)\frac{(x-\mu)^2}{2}) 
&=& 
\exp(-\frac{x^2}{2}) \exp(-(\tau+1)\frac{\mu^2}{2})
\exp((\tau+1) x\mu - \tau\frac{x^2}{2})  \nonumber \\
&=&
\exp(-\frac{x^2}{2}) \exp(-(\tau+1)\frac{\mu^2}{2})
\Big(1 + ((\tau+1) x\mu - \tau\frac{x^2}{2}) 
\nonumber \\ && \qquad\qquad 
+  \frac{1}{2}
((\tau+1) x\mu - \tau\frac{x^2}{2})^2 + \cdots \Big).
\label{eq:mixture-expansion}
\end{eqnarray}
The term $\exp(-(\tau+1)\frac{\mu^2}{2})$ can be absorbed into
$h(\mu,\tau;\theta)$ and can be ignored.  Also the constant term (i.e.\
terms not involving $x$) in the expansion can be ignored.
Now from (\ref{eq:convenient-condition-at-normal}) of Appendix
\ref{app:parameterization} it follows that without loss of generality we
can choose the prior distribution in such a way that the expected values
of the coefficients of $x$ and $x^2$ vanish.  
%
%
%
%
Therefore in (\ref{eq:mixture-expansion}) we only need to consider 
the cubic or higher degree terms in $x$ in the expansion.  
Relevant terms on the right-hand side of
(\ref{eq:mixture-expansion}) 
are
\begin{equation}
\label{eq;relevant-term}
\exp(-\frac{x^2}{2})\big[
-\frac{1}{2}\mu\tau x^3 + \frac{1}{8}\tau^2 x^4 + \frac{1}{6}\mu^3 x^3
-\frac{1}{4}\mu^2 \tau x^4 + \frac{1}{24}\mu^4 x^4
\big] . 
\end{equation}
If only the scale parameter is mixed, i.e.\ if $\mu\equiv 0$, then the
dominant term is $(1/8)\tau^2 x^4$. The primary example of this case 
is the family of $t$-distributions with $m=1/\theta$ degrees of
freedom, where the mixing distribution for the scale is the inverse
Gamma distribution. {}From the above consideration it follows that the
LBI test against the $t$-family is given by the standardized sample kurtosis.
On the other hand if only the location parameter is mixed, i.e.\ $\tau \equiv
0$ and $E_g(\mu^3)\neq 0$, then the LBI test is given by the
standardized sample skewness. 

More interesting case is that $\mu$ and $\tau$ is of the same order
and the LBI test involves both skewness and kurtosis simultaneously.  This 
happens in a limiting case of ``normal variance-mean mixture.''  In the
normal variance-mean mixture,  $X$ given $Y=y$ is normal with mean 
$a + b y$, $b\neq 0$, and variance $y$:
\[
X\; |\; Y=y \ \sim \ N(a+ b y, y), \qquad Y\sim g(y,\theta).
\]
Now assume that $Y$ degenerates to a constant as $\theta \rightarrow 0$.
Since we are considering location-scale invariant tests, we can assume
that $Y \rightarrow 1$ in distribution  and $a=-b$.
Writing $\mu=b(y-1)$ we have 
\begin{equation}
\label{eq:mu-tau}
\tau = \frac{1}{y} - 1 = -\frac{\mu}{\mu +b}= -\frac{\mu}{b} +
o(|\mu|), \qquad \text{or} \quad \mu = -b \tau + o(|\tau|).
\end{equation}
Therefore $\mu$ and $\tau$ become proportional as $\theta \rightarrow 0$.
In the following subsection we look at the 
generalized hyperbolic distribution as an example of this case.

\subsubsection{The case of the generalized hyperbolic distribution}
Generalized hyperbolic distribution (GH distribution)
was introduced by \cite{Barn:expo:1977}.
Detailed explanations including applications of GH distributions are given in
\cite{Barn-Shep:2001}, \cite{Eberlein:2001} or \cite{masuda2002}.
{}From \cite{Eberlein:2001} the density is written as
\begin{align}
\label{eq:density-GH}
&f_{GH}(x;\lambda,\alpha,\beta,\delta,\mu) \nonumber \\
&\qquad  =a(\lambda,\alpha,\beta,\delta)
\left(\delta^2+(x-\mu)^2\right)^{(\lambda-\frac{1}{2})/2}
K_{\lambda-\frac{1}{2}}\left(\alpha\sqrt{\delta^2+(x-\mu)^2}\right)\exp\left(\beta(x-\mu)\right),
\end{align}
where 
\[ a(\lambda,\alpha,\beta,\delta) =
\frac{(\alpha^2-\beta^2)^{\lambda/2}}{\sqrt{2\pi}
  \alpha^{\lambda-\frac{1}{2}}\delta^\lambda
  K_\lambda(\delta\sqrt{\alpha^2-\beta^2})}
\]
is the normalizing
constant and $K_\lambda$ is the modified Bessel function of the third
kind with index $\lambda$:
\[
K_\lambda(z)=\frac{1}{2}\int_0^\infty y^{\lambda-1}
\exp\left(-\frac{1}{2}z\left(y+y^{-1}\right)\right)dy,  \quad z > 0.
\]
The parameter space is given by 
\[
-\infty< \mu, \lambda < \infty, \quad \alpha > | \beta |, 
\]
with the additional boundaries $\{\delta=0, \lambda > 0\}$ and
$\{\alpha=|\beta|, \lambda < 0\}$.

GH distribution can
be characterized as a normal variance-mean mixture using the generalized
inverse Gaussian distributions (GIG distributions) as the mixing
distribution.  Let $X\ |\ Y=y$ be distributed as $N(\mu +\beta y, y)$
and let $Y$ have the generalized inverse Gaussian
distribution with parameters $\lambda$, $\delta$, and 
$\gamma=\sqrt{\alpha^2 - \beta^2}$. 
The density of $Y$ is written as 
\begin{equation}
\label{eq:density-GIG}
f_{GIG}(y;\lambda,\delta,\gamma)=\left(\frac{\gamma}{\delta}\right)^\lambda
\frac{1}{2K_\lambda(\delta\gamma)}y^{\lambda-1}
\exp\left(-\frac{1}{2}\left(\frac{\delta^2}{y}+\gamma^2 y\right)\right),
\quad y > 0,
\end{equation}
where the parameter space is given by
$
\gamma, \delta > 0, \ -\infty < \lambda < \infty, 
$
with the additional boundaries $\{\delta=0, \lambda > 0\}$ and
$\{\gamma=0,\lambda<0\}$.

In  (\ref{eq:density-GIG}) let
$\delta \to \infty$ and $\gamma \to \infty$ such that $\gamma/\delta \to
\overline{c}$, then 
it is easily seen that $Y$ degenerates to $\bar c$.
Therefore GH distribution converges to 
$N(\mu+\beta \overline{c},\overline{c})$ 
as $\delta \to \infty$ and $\gamma \to \infty$ such that $\gamma/\delta \to
\overline{c}$. 
As above we can assume
$\bar c=1$ and $\mu =-\beta$ without loss of generality.  
We also assume that $\beta$ is fixed. For
simplicity let $\delta=\gamma$. 
Then
(\ref{eq:density-GIG}) is written as
\[
f_{GIG}(y;\lambda,\gamma)=
\frac{1}{2K_\lambda(\gamma^2)}y^{\lambda-1}
\exp\left(-\frac{\gamma^2}{2}\left(\frac{1}{y}+y\right)\right).
\]
Note that this density has exponentially small tails at $y=0$ and
$y=\infty$.  Therefore term by term integration 
in (\ref{eq:mixture-expansion}) is justified.   

By (\ref{eq:mu-tau}), the main term in 
(\ref{eq;relevant-term}) is simply given as
\[
\exp(-\frac{x^2}{2})\big(
\frac{\beta}{2} x^3 + \frac{1}{8} x^4\big) \tau^2.
\]
It follows that the rejection region of the LBI test (for a fixed
$\beta$) is given by
\begin{eqnarray*}
c_{n+2}\sum_{i=1}^{n}z_i^4+4 \beta c_{n+1}\sum_{i=1}^{n} z_i^3 > k,
\end{eqnarray*}
where
\begin{equation}
\label{eq:gamma-l}
c_l = \int_0^\infty x^l e^{-nx^2/2} dx = \frac{2^{(l-1)/2}}{n^{(l+1)/2}}
\Gamma\big(\frac{l+1}{2}\big).
\end{equation}
We see that the LBI test involves both the skewness and the kurtosis
simultaneously and the weight depends on the value of $\beta$.

\subsection{Infinitely divisible family}
\label{subsec:id}

Here we consider an infinitely divisible family  with finite variance.
The characteristic function of an infinitely divisible random variable
$X$ with
mean 0 and variance 1 can be written as
\begin{equation}
\label{eq:id-ch}
\phi(t)=\exp [\int_{-\infty}^\infty (e^{itu} - 1 - itu) \frac{1}{u^2}
\mu (du)], 
\end{equation}
where the L\'evy measure $\mu$ can be taken as a probability measure.
Here we assume that $X$ possesses moments up to an appropriate order.
Since moments of the  L\'evy measure $\mu$  are the cumulants of $X$, 
existence of moments of $X$ up to an appropriate order is equivalent
to the existence of moments of $\mu$ to the same order.
For example if $Y$ has the exponential distribution, the
characteristic function of $X=Y-1$ can be written as
(\ref{eq:id-ch}) with $\mu(du)=u e^{-u}, u>0,$ (Example 8.10 of 
\cite{KSato-levy-processes}) and for the 
double-exponential distribution with variance 1, 
$\mu(du)=|u| e^{-\sqrt{2}u}, -\infty < u < \infty$.

Now we introduce the time parameter
$m=1/\theta$ and consider a L\'evy process  $X(m)$, where $X=X(1)$
has the characteristic function (\ref{eq:id-ch}).
Furthermore we standardize the variance as $X(m)/\sqrt{m}$.  Then by
the central limit theorem $X(m)/\sqrt{m}$ converges to $N(0,1)$ as
$m\rightarrow\infty$.  The characteristic function of 
$X(m)/\sqrt{m}$ is written as
\begin{equation}
\label{eq:id-ch-m}
\phi_m(t)=\phi(t/\sqrt{m})^m = \exp[\int_{-\infty}^\infty 
m (e^{iut/\sqrt{m}} - 1 - \frac{iut}{\sqrt{m}}) \frac{1}{u^2}
\mu (du)].
\end{equation}
Recalling the fact 
$
|e^{ix} - (1 + ix + (ix)^2/2 + \cdots + (ix)^k)/k|
\le |x|^{k+1}/(k+1)! 
$ for all real $x$, 
we can expand the integrand in (\ref{eq:id-ch-m}) as
\[
m (e^{iut/\sqrt{m}} - 1 - \frac{iut}{\sqrt{m}}) = -\frac{t^2}{2} + 
 \frac{(it)^3}{6 \sqrt{m}}u + \frac{(it)^4}{24m}u^2 + o(1/m)
\]
up to an appropriate order and integrate it term by term. Then
\begin{equation}
\label{eq:phim}
\phi_m(t)=\exp\Big( -\frac{t^2}{2} + \frac{\kappa_3}{6 \sqrt{m}} (it)^3 + 
  \frac{\kappa_4}{24m}(it)^4\Big) (1+ o(1/m)), 
\end{equation}
where $\kappa_j = \int_{-\infty}^\infty u^j \mu(du)$ is the $j$-th
cumulant of $X$. Note that (\ref{eq:phim}) is formally the same as the
usual Edgeworth expansion of the cumulant generating function of $m$
i.i.d.\ random variables.  By considering a L\'evy process, we can
allow $m$ to be fractional and we have a family of distributions $\{
X(m)/\sqrt{m}\}$ indexed by the continuous parameter $m=1/\theta$. By
the usual Edgeworth expansion, the density function of $X(m)/\sqrt{m}$ is
given as
\[
f(x;1/m) = \frac{1}{\sqrt{2\pi}} e^{-x^2/2} ( 1 + \frac{\kappa_3}{6\sqrt{m}}
H_3(x) + \frac{\kappa_4}{24m} H_4(x) + \frac{\kappa_3^2}{72m} H_6(x)) + o(1/m),
\]
where 
$H_j(x)$  
is the $j$-th Hermite polynomial.  We now see that i) if $\kappa_3
\neq 0$ then the LBI test is given by the sample skewness and ii) if
$\kappa_3=0$ and $\kappa_4 \neq 0$ then the LBI test is given by the
standardized sample kurtosis.

As examples consider the  centered exponential distribution and
the double-exponential distribution discussed at the beginning of this
section.  In the former case we test normality against the family of
normalized Gamma distributions and the LBI test is given by the
standardized sample skewness.  In the latter case, the characteristic
function of $X(m)/\sqrt{m}$ is given by
\[
\phi_m(t)=\Big(1-\frac{t^2}{2m}\Big)^{-m}
\]
This is a dual family of distributions to $t$-family in the sense of 
\cite{dreier-kotz-2002}. The LBI test against this family is given by
the sample kurtosis, as in the case of $t$-family.


\section{Testing against the stable family}
\label{sec:stable}

In this section as an important non-regular case we consider testing
against the stable family. 
The characteristic function of a general stable distribution
($\alpha \neq 1$) is  given by 
\begin{equation*}
\label{eq:characteristic-function}
\Phi(t)=\Phi(t;\mu,\sigma,\alpha,\beta)=
\exp\left(-|\sigma t|^{\alpha}\left\{1 + i\beta(\mbox{\rm sgn}
t)\tan\left(\frac{\pi\alpha}{2}\right)(|\sigma t|^{1-\alpha}-1)\right\} + i\mu t\right),
\end{equation*}
where $\mu$ is the location, $\sigma$ is the scale, $\beta$ is 
the ``skewness'' and $\alpha$ is the characteristic exponent. The 
parameter space is given by
\begin{equation*}
-\infty < \mu < \infty,\ \sigma >0,\ 0<\alpha \le2,\ |\beta|\le 1.
\end{equation*}
For the standard case
$(\mu,\sigma)=(0,1)$ we simply write the characteristic function as
\begin{equation}
\label{eq:characteristic-function-standard}
\Phi(t;\alpha,\beta)=\exp\left(-|t|^{\alpha}\left\{1 + i\beta(\mbox{\rm 
sgn}
t)\tan\left(\frac{\pi\alpha}{2}\right)(|t|^{1-\alpha}-1)\right\} \right).
\end{equation}
This is Zolotarev's (M) parameterization (see p.11 of \cite{Zolotarev:1986}).
The corresponding density is written as
$g(x;\mu,\sigma,\alpha,\beta)$  and $g(x;\alpha,\beta)$ in the standard
case. 

Letting $\alpha=2$ in
(\ref{eq:characteristic-function-standard}) we obtain $N(0,2)$. For
convenience let $\theta=2-\alpha$, $\mu=a$, $\sigma=b$ 
and we write \[
f(x;\theta)=g(x;a,b,2-\theta,\beta), 
\]
where $f(x;0)$ corresponds to $N(0,2)$.  For this section we take
$N(0,2)$ as the standard member of the normal location-scale family.
In the following we fix $\beta$ and for each $\beta$ we consider LBI
for $H_0:\theta=0$ vs\ $H_1:\theta > 0$.  This is similar to the case
of generalized hyperbolic distributions.  In particular for $\beta=0$
we are testing normality against the symmetric stable family, which is
important in practice.

It can be shown that we can differentiate
$
g(x;\alpha,\beta)= \frac{1}{2\pi} \int_{-\infty}^\infty e^{-itx}
\Phi(t;\alpha,\beta)dt
$
under the integral sign and the score function is written as
\begin{equation}
\label{eq:stable-score}
l_\theta(x;0)= -\frac{1}{2\pi} \int_{-\infty}^\infty e^{-itx}
\frac{\partial}{\partial \alpha}
\Phi(t;2,\beta)dt.
\end{equation}
In particular for $\beta=0$
\[
l_\theta(x;0)= \frac{1}{2\pi} \int_{-\infty}^\infty \cos(tx)
\log |t| \; t^2 e^{-t^2} dt .
\]

The non-regularity of stable family lies in the fact that this score
function has a very heavy tail.  In fact in \cite{Matsui:2005} it is
shown
that for large $|x|$
\[
l_\theta(x;0) = O \Big(\exp\big(\frac{x^2}{4}\big) |x|^{-3} 
\Big).
\]
Thus under $N(0,2)$, $E[l_\theta(x;0)]=0$  exists but
$E[l_\theta(x;0)^2]=\infty$ diverges.  This corresponds to
the fact that as $\alpha \uparrow 2$, the Fisher information
$I_{\alpha\alpha}$ diverges to infinity.  
\cite{Matsui:2005} gives a detailed analysis of the Fisher information
matrix for the general stable distribution close to the normal
distribution.


Although Assumption \ref{ass:1} does not hold for this case and we
have to give a separate proof, the following theorem holds.

\begin{theorem} \quad  
\label{thm:stable}
In the general stable family consider testing $H_0: \alpha=2$ vs.\  $H_1: \alpha<2$  
for fixed $\beta$.
Then the locally best invariant
is given by  (\ref{eq:lbi}), where the score function is given in 
(\ref{eq:stable-score}).
\end{theorem}

The proof of this theorem is very technical and it is given in
Appendix \ref{app:proof2}.  Note that score function puts extremely
heavy weights to outlying observations and this test can be considered
as an outlier detection test.  This is intuitively reasonable, because
the stable distribution with $\alpha < 2$ does not possess a finite
variance.

\section{Tests based on the profile likelihood}
\label{sec:profile}

In this section we consider tests based on the profile likelihood, where
the location and the scale parameters are estimated by the maximum likelihood.
We show that the LBI test and the test based on the profile likelihood are
different in general except for the case that the score function is a
third degree polynomial or a fourth degree polynomial without odd degree
terms.  Our argument in this section is formal and we implicitly
assume enough regularity conditions so that our formal argument is
justified.

Consider a density close to a normal distribution of the form
\begin{equation}
\label{eq:profile-likelihood-density}
f(x;\theta)= \frac{1}{\sqrt{2\pi}\sigma} \exp\Big(-
\frac{(x-\mu)^2}{2\sigma^2} \Big) \big\{ 1 + \theta
h\big(\frac{x-\mu}{\sigma}\big) + o(\theta)\},
\end{equation}
where $h$ is some smooth function. 

We estimate $\mu$ and $\sigma$ by the maximum likelihood under the null
and under the alternative and take the ratio of the maximized likelihoods.
$\theta$ is considered to be fixed in the estimation.  Since the
maximum likelihood estimator is location-scale equivariant, we obtain
an invariant testing procedure. Under the null hypothesis of normal
distribution the maximum likelihood estimates are $\hat \mu = \bar x$
and $\hat \sigma^2 = s^2$.  
Under the alternative, an approximation to $\hat \mu$
and $\hat \sigma^2$ to the order of $O(\theta)$ is easily derived as
\begin{equation}
\hat \mu = \bar x - \theta  s \frac{1}{n} \sum_{i=1}^n h'(z_i) 
+ o(\theta),\qquad 
\hat\sigma^2 = s^2 ( 1 - \theta \frac{1}{n} \sum_{i=1}^n  z_i
h'( z_i)) + o(\theta).
\label{eq:mle-theta}
\end{equation}
Let $L(\hat \mu, \hat \sigma^2)$ denote the log-likelihood under the
alternative and let $L(\bar x, s^2)$ denote the log-likelihood under
the null.  Then  substituting (\ref{eq:mle-theta}) into
(\ref{eq:profile-likelihood-density}) we obtain
\begin{align*}
L(\hat\mu, \hat\sigma^2)
&= -\frac{1}{2}\sum_{i=1}^n (x_i-\bar x)^2 \frac{1}{s^2}
( 1 + \theta \frac{1}{n} \sum_{i=1}^n  z_i ) - n \log s
+ \frac{\theta}{2} \frac{1}{n}\sum_{i=1}^n  z_i h'( z_i) + 
o(\theta)\\
&= L(\bar x, s^2) - \frac{\theta}{2}\frac{n-1}{n} \sum_{i=1}^n 
z_i h'( z_i) + o(\theta).
\end{align*}
Hence the test based on the profile likelihood ratio has the rejection
region
\begin{equation}
\label{eq:profile}
\sum_{i=1}^n 
z_i h'( z_i) > k .
\end{equation}
On the other hand, as discussed in Section \ref{sec:lbi},
for large $n$ the Laplace approximation to the integral in 
(\ref{eq:lbi}) implies that the LBI is asymptotically equivalent to 
\begin{equation}
\label{eq:lbi-asymp}
\sum_{i=1}^n h( z_i) > k' .
\end{equation}
We see that (\ref{eq:profile}) and (\ref{eq:lbi-asymp}) are generally
different even asymptotically.  It should be noted that if 
$h$ is a third degree polynomial or a fourth degree polynomial without
odd degree terms, then both the profile likelihood procedure and the LBI
procedure reduce to the sample skewness and the sample kurtosis.

\section{Multivariate extensions}
\label{sec:multivariate}

In this section we consider multivariate extensions of our results.  A
comprehensive survey on invariant tests of multivariate normality is
given in \cite{Henze:2002}.


For a column vector $a\in \R^p$ and a $p\times p$ nonsingular matrix $B$,
let
\begin{equation}
\label{multivariate:f}
 f_{a,B}(x;\theta) = \frac{1}{|\det B|} f(B^{-1}(x-a);\theta), \quad
 x\in \R^p, 
\end{equation}
be a one-parameter family with the shape parameter $\theta$.
As in the univariate case, we assume that
$$ f(x;0) = \frac{1}{(2\pi)^{p/2}} \exp(-\Vert x\Vert^2/2), $$
where $\Vert\cdot\Vert$ denotes the standard Euclidean norm in $\R^p$.
Based on the i.i.d.\ samples $x_1,\ldots,x_n$ from $f_{a,B}(x;\theta)$,
we discuss invariant testing procedures for testing the normality
$$
H_0: \theta=0 \quad {\rm vs.}\quad  H_1: \theta>0.
$$

Write $X=(x_1,\ldots,x_n)'\in \R^{n\times p}$.
Consider the group
$$ \R^p\times GL(p) 
= \{ (a,B) \mid a\in \R^p,\ B\in\R^{p\times p},\ \det B \ne 0 \} $$
endowed with the product 
$(a_1,B_1) \cdot (a_2,B_2) = (B_2 a_1 + a_2, B_2 B_1)$.
This group acts on the sample space $\R^{n\times p}$ as
\begin{equation}
\label{multivariate:action}
 (a,B) X = 1_n a' + XB', \quad (a,B)\in \R^p\times GL(p),
\end{equation}
where $1_n=(1,\ldots,1)'\in \R^n$.
For each $\theta$ fixed, the action (\ref{multivariate:action}) 
induces the transitive action on the parameter space.
In other words, the model (\ref{multivariate:f}) is a transformation 
model with the parameter $(a,B)$.
Thus, it is natural to consider invariant procedures under the action
 (\ref{multivariate:action}).   

Let $LT(p)$ be the set of $p\times p$ lower triangular matrices with
positive diagonal elements.
Let $\bar x=\sum_{i=1}^n x_i/n$ and
 $S=\sum_{i=1}^n (x_i-\bar x) (x_i-\bar x)'/n$ be
the sample mean vector and the sample covariance matrix.
Let $T\in LT(p)$ be the Cholesky root of $S$ so that $S=TT'$.
Let
\begin{equation}
\label{multivariate:invariant}
 Z = (z_1,\ldots,z_n)' = (X - 1_n\bar x')(T')^{-1}
\end{equation}
($z_i=T^{-1}(x_i-\bar x)$, $i=1,\ldots,n$).
It is easy to see that a maximal invariant
 under the action (\ref{multivariate:action}) is
$$ W = ZZ' =  (X - 1_n\bar x') S^{-1} (X - 1_n\bar x')', $$
and we can choose a cross section
$\tilde Z=\tilde Z(X)=(\tilde z_1,\ldots,\tilde z_n)' \in \R^{n\times p}$
as a unique decomposition of $W=\tilde Z\tilde Z'$ in some appropriate way.
Note that $\tilde Z = ZQ'$, or $\tilde z_i=Q z_i$, for some
$p\times p$ orthogonal matrix $Q$.
The following is a multivariate extension of 
 Proposition \ref{eq:basic-prop}.

\begin{proposition}\quad
\label{multivariate:prop1}
Under the group action of $\R^p\times GL(p)$,
the critical region of the most powerful invariant test for testing
$H_0: \theta=0$ against $H_1: \theta = \theta_1 > 0$ is given by
\begin{equation}
\label{eq:lambdart}
\frac{\int_{GL(p)} \int_{\R^p} \prod_{i=1}^n f(a + B z_i ; \theta_1)
  |\det B|^{n-p-1} da  dB}
     {\int_{GL(p)} \int_{\R^p} \prod_{i=1}^n f(a + B z_i ; 0)
  |\det B|^{n-p-1} da  dB} > k
\end{equation}
for some $k > 0$, 
 where $da=\prod_{i=1}^p da_i$ and $dB=\prod_{i,j=1}^p db_{ij}$
are the Lebesgue measures of $\R^p$ and $\R^{p\times p}$, respectively.
\end{proposition}

\noindent
{\it Proof\/}.
The Jacobian of the transformation $X\mapsto (a,B)X = 1_n a' + X B'$
is $(\det B)^n$.
The left invariant measure of $\R^p\times GL(p)$ is $(\det B)^{-(p+1)} da dB$.
From Theorem 4 of \cite{Wijsman:1967}, the critical region is
$$
\frac{\int_{GL(p)} \int_{\R^p} \prod_{i=1}^n f(a + B \tilde z_i ; \theta_1)
  |\det B|^{n-p-1} da  dB}
     {\int_{GL(p)} \int_{\R^p} \prod_{i=1}^n f(a + B \tilde z_i ; 0)
  |\det B|^{n-p-1} da  dB} > k,
$$
which is equivalent to  (\ref{eq:lambdart}).
\hfill $\Box$
\bigskip

Next consider the subgroup
$$ \R^p\times LT(p) = \{ (a,T) \mid a\in \R^p,\ T \in LT(p) \} $$
of $\R^p\times GL(p)$.
This also acts on the sample space $\R^{n\times p}$ 
with the same action (\ref{multivariate:action})
with $GL(p)$ replaced by $LT(p)$.
For this group, the induced action on the parameter $(a,B)$ in the model
 (\ref{multivariate:f}) is not transitive anymore.
However, when we consider a subclass
of (\ref{multivariate:f}) that
\begin{eqnarray}
 f_{a,B}(x;\theta)
&=& \frac{1}{|\det B|} h(\Vert B^{-1}(x-a)\Vert^2 ;\theta)
\label{multivariate:h} \\
&=&  \frac{1}{\sqrt{\det (BB')}} h((x-a)'(BB')^{-1}(x-a) ;\theta)
\nonumber
\end{eqnarray}
($h$ is a function), the action on the parameter $(a,BB')$ is transitive,
and invariant testing procedures under the group $\R^p\times LT(p)$
may be more appropriate in some cases.

For the action of $\R^p \times LT(p)$, 
$Z$ in (\ref{multivariate:invariant}) is a maximal invariant, 
and we can use $Z$ itself as a cross section.

The most powerful invariant test under the action of $\R^p\times LT(p)$
is given as follows.
\begin{proposition}\quad
\label{multivariate:prop2}
Under the group action of $\R^p\times LT(p)$,
the critical region of the most powerful invariant test for testing
$H_0: \theta=0$ against $H_1: \theta = \theta_1 > 0$ is given by
$$
\frac{\int_{LT(p)} \int_{\R^p}
      \prod_{i=1}^n f(a + T z_i ; \theta_1) da
      \prod_{i=1}^p t_{ii}^{n-i-1} \prod_{i\ge j}dt_{ij}}
     {\int_{LT(p)} \int_{\R^p}
      \prod_{i=1}^n f(a + T z_i ; 0) da
      \prod_{i=1}^p t_{ii}^{n-i-1} \prod_{i\ge j}dt_{ij}} > k'
$$
for some $k' > 0$, where $T=(t_{ij})\in LT(p)$.
\end{proposition}

\noindent
{\it Proof\/}.
The Jacobian of the transformation $X\mapsto (a,T)X = 1_n a' + X T'$ is
 $(\det T)^n=\prod_{i=1}^n t_{ii}^n$.
The left invariant measure of $\R^p\times LT(p)$ is 
 $da \prod_{i=1}^p t_{ii}^{-(i+1)} \prod_{i\ge j}dt_{ij}$.
The proposition follows from Theorem 4 of \cite{Wijsman:1967}.
\hfill $\Box$
\bigskip

{}From Propositions \ref{multivariate:prop1} and \ref{multivariate:prop2},
under similar conditions to Assumption \ref{ass:1}, the LBI
test can be derived by integrating the score function
$\sum_{i=1}^n \ell_\theta (a+B z_i; 0)$ with respect to $(a,B)$.
In the rest of this section, we examine a particular case where
$$ \prod_{i=1}^n f(a+Bz_i;\theta) = \prod_{i=1}^n f(a+Bz_i;0)
 \{1 + \theta \sum_{i=1}^n \ell_\theta(a+Bz_i;0) + o(\theta) \} $$
with
$$ \ell_\theta(x;0) = p_0 \Vert x\Vert^4 +  p_1 \Vert x\Vert^2 +  p_2. $$
This holds, for example, when
$f_{a,B}(x;\theta)$ is of the form of (\ref{multivariate:h})
with
$$ h(y;\theta) = 
\begin{cases}
 \displaystyle
 \frac{\Gamma((p+\theta^{-1})/2)}
 {(\theta^{-1}\pi)^{p/2} \Gamma(\theta^{-1}/2)}
 (1+\theta y)^{-(p+\theta^{-1})/2} & (\theta>0) \\
 \displaystyle
 \frac{1}{(2\pi)^{p/2}} \exp(-y/2) & (\theta=0)
\end{cases}
$$
(multivariate $t$ distribution with $\theta^{-1}$ degrees of freedom).
We restrict our attention to the case $p_0>0$ for simplicity.
\begin{assumption} \quad
\label{ass:2}

\noindent
(i)
$$ 
  \frac{\frac{\partial}{\partial \theta} f(a+Bz ; \theta)|_{\theta=0}}
  {f(a+Bz ;0)} = p_0 \Vert z\Vert^4 +  p_1 \Vert z\Vert^2 +  p_2
 \quad (p_0>0). $$

\noindent
(ii)\ \ %
For some $\epsilon>0$, 
$$
\int_{GL(p)} \int_{\R^p} g(a,B;\epsilon)^n \exp(-\frac{n}{2}\Vert a\Vert^2
 - \frac{n}{2}\tr(B'B)) |\det B|^{n-p-1} da dB < \infty,
$$ 
where
$$
g(a,B;\epsilon)=\sup_{\Vert z\Vert\le 1,\,0\le \theta \le \epsilon} 
  \frac{ |\frac{\partial}{\partial \theta} f(a+Bz ;
    \theta)|}{f(a+Bz ;0 )}.
$$
\end{assumption}  
%
\begin{theorem}
\label{multivariate:thm}
Under Assumption \ref{ass:2}, the rejection region of the LBI
test for testing normality 
$H_0: \theta=0$ vs.\ $H_1: \theta=\theta_1>0$
under the action of $\R^p\times GL(p)$ is given by
$$ \sum_{i=1}^n \Vert z_i\Vert^4 >k. $$
The rejection region of the LBI test 
 under the action of $\R^p\times LT(p)$ is given by
\begin{eqnarray*}
&& (n+p+2)(n+p) \sum_{i=1}^n \Vert z_i\Vert^4
 - 2(n+p+2) \sum_{i=1}^n \sum_{j,k=1}^p \max(j,k) z_{ij}^2 z_{ik}^2 \\
&& \qquad
 - 2(n+p) \sum_{i=1}^n \sum_{j,k=1}^p \min(j,k) z_{ij}^2 z_{ik}^2
 + 4 \sum_{i=1}^n (\sum_{j=1}^p j z_{ij}^2)^2 > k',
\end{eqnarray*}
where $z_{ij}$ is the $j$th element of $z_i$.
\end{theorem}

The lemma below is used in proving Theorem \ref{multivariate:thm}.
This is easily proved by some standard Jacobian formulas
in the multivariate analysis (e.g., page 86 of \cite{Muirhead:1982}).

\begin{lemma}
\label{lem:wishart}
Let $Sym(p)$ denote the set of $p\times p$ real symmetric matrices.
Define a map $\varphi: \R^{p\times p} \to Sym(p)$ by
 $\varphi(B) = n B'B$.  Then, for any measurable set $D\subset Sym(p)$,
$$ \int_{\varphi(B)\in D}  \exp( -\frac{n}{2}\tr(B'B)) |\det B|^{n-p-1} dB
 \propto
   \int_D  \exp( -\frac{1}{2}\tr\,C) (\det C)^{\frac{1}{2}(n-p-2)} dC,
$$
where $C=(c_{ij})\in Sym(p)$ and $dC = \prod_{i\ge j} dc_{ij}$.
\end{lemma}

\noindent
{\it Proof of Theorem \ref{multivariate:thm}\/}.
Note first that
$\sum_{i=1}^n \Vert a + B z _i \Vert^2 = n \Vert a\Vert^2 + n \tr(B'B)$ 
because $\sum_{i=1}^n z_i=0$ and $\sum_{i=1}^n z_i z_i' = n I_p$.
The second and the third terms of $\ell_\theta$ are irrelevant to
$z_i$'s.

In the case of $\R^p\times GL(p)$, the rejection region is
of the form $\sum_{i=1}^n I(z_i) >k$, where
$$ I(z) = \int_{GL(p)} \int_{\R^p} \Vert a + B z \Vert^4
 \exp(-\frac{n}{2}\Vert a\Vert^2 -\frac{n}{2}\tr(B'B)) |\det B|^{n-p-1}
 da dB. $$
By Lemma \ref{lem:wishart}
the integral of a function of $B'B$ can be
replaced by taking expectation with respect to the Wishart distribution
$n B'B \sim W_p(n-1,I_p)$.  
On the other hand, the integration with respect
to $a$ is regarded as the expectation with respect to
 $\sqrt{n}a\sim N_p(0,I_p)$.
Note that for the Wishart matrix $C \sim W_p(n-1,I_p)$, it holds that
$$ E[z'C z] = (n-1) \Vert z\Vert^2, \qquad
 E[(z'C z)^2] = (n-1)(n+1) \Vert z\Vert^4. $$
By taking expectations of
\begin{eqnarray*}
\Vert a + B z \Vert^4
&=& (\Vert a\Vert^2 + 2 a'B z + z'B'B z)^2 \\
&=& (z'B'B z)^2 + 2 \Vert a\Vert^2 (z'B'B z) \\
&& + (\mbox{terms of odd degrees in $a$})
   + (\mbox{a term independent of $z$}),
\end{eqnarray*}
we see that $\sum_{i=1}^n I(z_i)$ is proportional to
$\sum_{i=1}^n \Vert z_i\Vert^4 + \mbox{const.}$

In the case of $\R^p\times LT(p)$, the rejection region is
of the form $\sum_{i=1}^n I(z_i) >k'$, where
$$ I(z) = \int_{LT(p)} \int_{\R^p} \Vert a + T z \Vert^4
 \exp(-\frac{n}{2}\Vert a\Vert^2 -\frac{n}{2}\tr(T'T)) da
 \prod_{i=1}^n t_{ii}^{n-i-1} \prod_{i\ge j} dt_{ij}. $$
The integration with respect to $T$ is reduced to taking expectations
$n t_{ii}^2\sim\chi^2_{n-i-1}$ and $\sqrt{n}t_{ij}\sim N(0,1)$ ($i>j$).
The details are given in Appendix \ref{multivariate:appendix}.
\hfill $\Box$

\appendix
\section{Proofs}

\subsection{Proof of Theorem \ref{thm:1}}
\label{app:proof1}
By the mean value theorem
$$
f(x;\theta)=f(x;0) + \theta \frac{\partial}{\partial \theta}
f(x;\theta^*), 
$$
where $0 < \theta^* =\theta^*(x) < \theta$.  Then
\begin{eqnarray*}
\prod_{i=1}^n  f(a+b z_i;\theta)
&=&
\prod_{i=1}^n \left(f(a+b z_i;0) + \theta \frac{\partial}{\partial \theta}
f(a+bz_i;\theta^*)\right) \\
&=& \prod_{i=1}^n f(a+b z_i;0) \times (1 + \theta \sum_{i=1}^n 
\frac{\frac{\partial}{\partial \theta}
f(a+bz_i;\theta^*)}{f(a+bz_i;0)} + \theta^2 R)
\end{eqnarray*}
where 
$$
R = \sum_{l=2}^n \theta^{l-2}
\sum_{1\le i_1 < \cdots < i_l \le n} 
  \frac{ \frac{\partial}{\partial \theta} f(a+b z_{i_1} ;
    \theta^*)}{f(a+bz_{i_1} ;0 )} \cdots
  \frac{ \frac{\partial}{\partial \theta} f(a+b z_{i_l} ;
    \theta^*)}{f(a+bz_{i_l} ;0 )}.
$$
By Assumption \ref{ass:1}
\begin{equation}
\label{eq:remainder}
\int_0^\infty \int_{-\infty}^\infty |R| 
\exp\Big(-\frac{n(a^2+b^2)}{2}\Big) b^{n-2} da db < \infty.
\end{equation}
Furthermore by the continuous differentiability of $f(x;\theta)$ with
respect to $\theta$ and the dominated convergence theorem we have
\begin{eqnarray*}
&&
\int_0^\infty \int_{-\infty}^\infty 
\sum_{i=1}^n 
\frac{\frac{\partial}{\partial \theta}
f(a+bz_i;\theta^*)}{f(a+bz_i;0)}
\exp\Big(-\frac{n(a^2+b^2)}{2} \Big) b^{n-2} da db \\
&& \qquad \rightarrow
\int_0^\infty \int_{-\infty}^\infty 
\sum_{i=1}^n 
l_{\theta}(a+z_i;0)
\exp\Big(-\frac{n(a^2+b^2)}{2} \Big) 
b^{n-2} da db   \qquad (\theta \rightarrow 0).
\end{eqnarray*}
Now the theorem follows by the standard argument on the locally most
powerful test (e.g.\ Section 4.8 of \cite{Cox:Hink:1974}). \hfill $\Box$

\subsection{Proof of Theorem \ref{thm:stable}}
\label{app:proof2}

In the proof, $M>0$ denotes some suitable constant.
Since Assumption 1 is not applicable for Theorem 
\ref{thm:stable}, we have to prove  the finiteness
of  (\ref{eq:remainder}) by a separate argument.
It suffices to prove that for each subsequence
$1\le i_1 < \dots < i_l \le n$
\[
\int_0^\infty \int_{-\infty}^{\infty} 
 \left|\frac{\partial}{\partial \theta} f(a+b z_{i_1} ; \theta^*)\cdots
 \frac{\partial}{\partial \theta} f(a+b z_{i_l} ; \theta^*) \right |
\exp\left(-\frac{1}{4}\sum_{k\neq i_j}(a+bz_k)^2\right)b^{n-2} da db < \infty.
\]
Without loss of generality consider  $i_1 = 1,\dots, i_l=l$ 
and  write
\begin{equation*}
W_l=\left| \frac{\partial}{\partial \theta} f(a+b z_1 ; \theta^*) 
\cdots \frac{\partial}{\partial \theta} f(a+b z_l ; \theta^*)
\right|\exp\left(-\frac{1}{4}\sum_{k=l+1}^n(a+bz_k)^2\right)b^{n-2}.
\end{equation*}
For evaluations of $W_l$ we need the following property of the score
function of general stable distributions. It follows from Lemma 3.1
of \cite{Matsui:2005}.
\begin{lemma}
\label{lem:boundness-continuity-tail-order-derivative}
For $\alpha=2-\theta \neq 1$, 
$|(\partial/\partial \theta) f(x ; \theta) |$ is bounded and uniformly
continuous in $x$. Furthermore as $\theta=2-\alpha \downarrow 0$,
there exist $M>0$, $x_0>0$, such that
\[
\left|\frac{\partial}{\partial \theta} f(x ; \theta) \right|\le
 M\cdot|x|^{\theta-3} \log|x| ,
\qquad  \forall |x| \ge x_0.
\]
\end{lemma}
The integrability of $W_l$ for $l\le n-1$ follows from that 
of $W_n$, since  $\exp(-1/4x^2) \le
M\cdot |\partial/\partial \theta f(x;\theta^\ast)|$ 
from Lemma \ref{lem:boundness-continuity-tail-order-derivative}.
However, the integrability of $W_n$ needs a very detailed argument.
We replace $a$
by $r=a+bz_1, $ 
then $W_n$ 
becomes
\begin{equation}
\label{eq:g(r,b)}
W_n(r,b) \equiv
\prod_{k=1}^{n} 
\left|\frac{\partial}{\partial \theta} f(r+b(z_k-z_1) ;\theta^*)
\right|b^{n-2}.
\end{equation}
Note that $z_k-z_1\neq 0$ implies
$$
\exists c>0 \quad \mbox{s.t.} \quad \forall k\neq 1\quad |c(z_k-z_1)|>2,  
$$
\begin{equation}
\label{eq:tirianglar-inequality-x-large}
b>cx>c|r|\ \Rightarrow\ |r+b(z_j-z_1)| > x.
\end{equation}
Now we divide the integral of (\ref{eq:g(r,b)}) into three parts
\[
\left(\int_{|r| \le x_0}\int_0^\infty 
\quad +\int_{|r|>x_0}\int_{b\le c|r|} 
\quad +\int_{|r|>x_0}\int_{b> c|r|}\ \right) W_n(r,b) dr db 
\equiv   I_1+I_2+I_3.
\]

Using Lemma \ref{lem:boundness-continuity-tail-order-derivative} and (\ref{eq:tirianglar-inequality-x-large})
in $I_1$, we have
\begin{eqnarray*}
I_1 
&\le& \int_{|r| \le x_0}\int_{b\le c x_0} W_n(r,b) dr db\\
&& \qquad\qquad 
+M\cdot\int_{|r| \le x_0}\int_{b> c x_0} \max_k 
\left( \frac{ \log |r+b(z_k-z_1)| }{|r+b(z_k-z_1)|^{3-\theta^\ast}} \right)^{n-1}
b^{n-2} dr db < \infty
\\
&<& \infty.
\end{eqnarray*}

For $I_2$ the following lemma is useful.
\begin{lemma}
\label{lem:integrand-I2}
Suppose that $\{z_k\neq 0:k \in n,\ z_k \neq z_j\}$ are given.
Then 
\begin{equation}
\label{eq:lem-I2}
\frac{W_n(r,b)}{\left|
\frac{\partial}{\partial \theta} f(r ; \theta^*)\right|}  =
\prod_{k=2}^{n}
\left|\frac{\partial}{\partial \theta} f(r+b(z_k-z_1) ;\theta^*)
\right| \cdot
b^{n-2}
\end{equation}
is bounded in $-\infty< r <\infty $ and $b>0$.
\end{lemma}
\begin{proof}
Assume that  (\ref{eq:lem-I2}) is not bounded. Choose a sequence of
$(r,b)$ such that (\ref{eq:lem-I2}) diverges to $\infty$.
Since the terms in the absolute value on the right-hand side are 
bounded, $b$ has to go to $\infty$.
By the assumption we can choose $r$ such that for some $k$, 
$$
|r + b(z_k-z_1)| < c' b^{\gamma},
$$
where $c'>0$ is a constant and $0\le \gamma <1$ (otherwise 
(\ref{eq:lem-I2}) converges to $0$ as $b \uparrow\infty$ from Lemma \ref{lem:boundness-continuity-tail-order-derivative}).
Then for $k\neq l$ we
have
$$
|r+b(z_k-z_1)-\{r+b(z_l-z_1)\}| = b|z_k-z_l|.
$$
Hence as $b \uparrow \infty$ 
\begin{equation}
\label{eq:inequality-x-large}
|r+b(z_l-z_1)|>x_0.
\end{equation}
Furthermore, since
$c' b^\gamma <b|z_l-z_1|$ for sufficiently large $b$, the triangular equality gives
$$
 \left|\frac{r}{\sqrt{b}}+\sqrt{b}(z_l-z_1) \right|\ge 
\left|c' b^{\gamma-1/2}-\sqrt{b}|z_k-z_l| \right| \uparrow\ \infty,\quad
\mbox{as}\ b \uparrow \infty.
$$ 
Therefore, from Lemma
 \ref{lem:boundness-continuity-tail-order-derivative} and (\ref{eq:inequality-x-large}), as $b \uparrow \infty$ the left-hand side of (\ref{eq:lem-I2})
approaches 
\begin{eqnarray*}
&& \left| \frac{\partial}{\partial \theta} f(r+b(z_k-z_1) ; \theta^*)
\right|\cdot 
\prod_{l\neq k} \frac{\log|r+b(z_l-z_1)| }{|r+b(z_l-z_1)|^{3-\theta^*}}\cdot b^{n-2} \\
&&\le M\cdot \prod_{l\neq k}
 \frac{\log|r+b(z_l-z_1)|}{|r+b(z_l-z_1)|^{1-\theta^*}}\frac{1}{|r/\sqrt{b}+\sqrt{b}(z_l-z_1)|^2}\ \  \downarrow\  0,
\end{eqnarray*}
regardless of selection of $k$.
This is a contradiction and the proof is over.
\end{proof}
By Lemma \ref{lem:integrand-I2} we get
\begin{eqnarray*}
I_2 &&\le
\sup_{r,b} \left\{ W_n(r,b) \Big/ \left|
\frac{\partial}{\partial \theta} f(r ; \theta^*)\right|
\right\} \cdot 
\int_{|r|>x_0}\int_{b \le c|r|} \left|
\frac{\partial}{\partial \theta} f(r ; \theta^*)\right| db  dr \\
 &&\le M\cdot 
\int_{|r|>x_0} \left|
\frac{\partial}{\partial \theta} f(r ; \theta^*)\right|\cdot 2c|r|  dr 
\ < \infty.
\end{eqnarray*}

Finally for $I_3$ from Lemma
\ref{lem:boundness-continuity-tail-order-derivative} and (\ref{eq:tirianglar-inequality-x-large}), 
\begin{eqnarray*}
I_3 \le M\cdot \int_{|r|>x_0} \frac{\log{r}}{|r|^{3-\theta^*}} \int_{b>c|r|} \max_k
\left(\frac{\log |r+b(z_k-z_1)|}{|r+b(z_k-z_1)|^{3-\theta^* } }\right)^{n-1}
b^{n-2} db dr.
\end{eqnarray*}
Since for large $x>0$, $(\log x)^{n-1} \le x$, we have
$$
\int_{b>c|r|}
\left(\frac{\log |r+b(z_k-z_1)|}{|r+b(z_k-z_1)|^{3-\theta^* } }\right)^{n-1}
b^{n-2}db \le M\cdot \int_{b>c|r|}|r+b(z_k-z_1)|^{-(n-1)(3-\theta^*)+1} b^{n-2} db.
$$
The right-hand side 
is bounded by
the equation 2.111, 2 on p.67 of \cite{Gradshteyn:Ryzhik:2000}:
$$
\int \frac{x^l}{z_1^m}dx = \frac{x^l}{z_1^{m-1}(l+1-m)b'}-\frac{na'}{(l+1-m)b'}\int\frac{x^{l-1}}{z_1^m}dx,
$$ 
where $z_1=a'+b'x$ and $a'$, $b'$ are constants.
By induction we obtain 
\begin{eqnarray*}
\int \frac{x^l}{(a'+b'x)^m}dx &=&
 -\frac{x^l}{(m-l-1)(a'+b'x)^{m-1}b'} \\
&&-\sum_{k=1}^{l}
\frac{l(l-1)\cdots(l+1-k)a'^kx^{l-k}}{(m-l-1)\cdots(m-l-1+k)(a'+b'x)^{m-1}b'^{k+1}}.
\end{eqnarray*}
Letting $a'=r$, $b'=(z_k-z_1)$, $m=\lfloor(n-1)(3-\theta^*)\rfloor-1$, $l=n-2$,
in the equation above
and 
utilizing Lemma \ref{lem:boundness-continuity-tail-order-derivative}, 
we have
\begin{equation*}
\int_{b>c|r|}
\left(\frac{\log |r+b(z_k-z_1)|}{|r+b(z_k-z_1)|^{3-\theta^* } }\right)^{n-1}
b^{n-2}db \le M \cdot \frac{n-1}{|r|^{\lfloor(n-1)(3-\theta)\rfloor-n}}.
\end{equation*}
Since the right-hand side 
is integrable with respect
to $r>x_0$, we have $I_3 < \infty$.  This completes the proof.

\section{Question of parametrization}
\label{app:parameterization}

Here we briefly discuss how to choose $a(\theta)$ and $b(\theta)$ in
(\ref{eq:another-choice}).   
Write $l(x;\theta)=\log f(x;\theta)$.
Under the assumption that the $3\times 3$ Fisher
information matrix exists at $(a(\theta),b(\theta),\theta)$,
it is convenient to determine $(a'(\theta), b'(\theta))$ in such a way
that $(d/d\theta) l_{a(\theta),b(\theta)}(x;\theta)$ is orthogonal 
to the location-scale family in the sense of Fisher information, i.e.\ 
\begin{eqnarray}
\label{eq:diff-eq-1}
&& 
\int \frac{d}{d\theta} l_{a(\theta),b(\theta)}(x;\theta)  
\frac{\partial}{\partial a} l_{a(\theta),b(\theta)}(x;\theta)
f_{a(\theta),b(\theta)}(x;\theta) dx=0 \\
&&
\label{eq:diff-eq-2}
\int \frac{d}{d\theta} l_{a(\theta),b(\theta)}(x;\theta)  
\frac{\partial}{\partial b} l_{a(\theta),b(\theta)}(x;\theta)
f_{a(\theta),b(\theta)}(x;\theta) dx=0 
\end{eqnarray}
These give  a system of differential equations for
$(a(\theta),b(\theta))$.

Actually we are only concerned in the neighborhood of the normal
distribution and we only consider determining $(a'(0), b'(0))$.
At $\theta=0$, $l(x;0)= -(1/2)\log(2\pi) - x^2/2$.  Therefore
$$
\frac{\partial}{\partial a}
l_{a(\theta),b(\theta)}(x;\theta)_{|\theta=0}   = x, \quad
\frac{\partial}{\partial b}
l_{a(\theta),b(\theta)}(x;\theta)_{|\theta=0}    = x^2
$$
$$
\frac{d}{d\theta}l_{0,1}(x;\theta)=-\frac{1}{b'(0)} + b'(0) x^2 +
a'(0) + l_\theta(x;0)
$$
and 
(\ref{eq:diff-eq-1}), (\ref{eq:diff-eq-2}) reduce to
$$
\int (-\frac{1}{b'(0)} + b'(0) x^2 +a'(0)x + l_\theta(x;0)) 
 x^k \phi(x)dx=0, \qquad k=1,2,
$$
which can be solved for $a'(0)$ and $b'(0)$.  

Note that we do not necessarily have to explicitly solve for $a'(0)$
and $b'(0)$.  Instead for theoretical developments we can use the fact
that the standard member $f(x;\theta)$ can be chosen in such a way
that
\begin{equation}
\label{eq:convenient-condition-at-normal}
\int l_\theta(x;0) x^k \phi(x)dx=0, \qquad k=1,2.
\end{equation}
When $l_\theta(x;0)$ is a polynomial in $x$, 
(\ref{eq:convenient-condition-at-normal}) shows that we can choose 
$l_\theta(x;0)$ such that it is cubic or of higher degree in $x$.
This is enough for simplifying our treatment of mixing distribution in
Section \ref{subsec:mixture}.  

\section{Details in the case of  polynomial  score function}
\label{app:hermite}

Here we write out coefficients of LBI in the case of polynomial score
function (cf.\ (\ref{eq:polynomial-score})). 
Suppose that $l_\theta(x;0)$ is given as
\[
l_\theta(x;0) = c_0 x^k + c_1 x^{k-1} + \dots + c_k
  = \sum_{j=0}^k c_{k-j} x^j.
\]
Then 
\[
\sum_{i=1}^n l_\theta(a+bz_i;0) 
= \sum_{i=1}^n \sum_{j=0}^k c_{k-j} (a+bz_i)^j 
= n\sum_{j=0}^k c_{k-j}  \sum_{l=0}^j\binom{j}{l} a^l b^{j-l} 
 \tilde m_{j-l}.
\]
Using (\ref{eq:gamma-l})
for even $l$ we have
\[
\int_0^\infty \int_{-\infty}^\infty a^l  b^{j-l} \exp(-\frac{n(a^2+b^2)}{2})
  b^{n-2} da db
= \frac{2^{(n+j-2)/2}}{n^{(n+j)/2}} \Gamma\big(\frac{l+1}{2}\big) \times 
 \Gamma\big(\frac{n+j-l-1}{2}\big).
\]
For odd $l$ the integral is zero.  Also we only consider $j-l \ge 3$.
Hence the LBI test statistic is
given as
\begin{equation}
\label{eq:LBI-explicit}
\sum_{j=3}^k c_{k-j}  \left(\frac{2}{n}\right)^{(n+j-2)/2}
\sum_{l=0, \ l:{\rm even}}^{j-3}
\binom{j}{l} 
\Gamma\big(\frac{l+1}{2}\big) \times 
 \Gamma\big(\frac{n+j-l-1}{2}\big) \tilde m_{j-l}.
\end{equation}


\section{Moments of $z'T'Tz$}
\label{multivariate:appendix}

Let $T=(t_{ij})\in LT(p)$ be a random matrix whose diagonal and lower
off-diagonal elements are independently distributed as
$t_{ii}\sim\chi_{m+p-i}$, $t_{ij}\sim N(0,1)$ ($i>j$), where $m>0$ is a
constant.
Let $z=(z_1,\ldots,z_p)'\in \R^p$ be a constant vector.
In this section we evaluate the expectations
$$ R_p(z) = E[z'T'Tz], \qquad S_p(z) = E[(z'T'Tz)^2] $$
required in proving Theorem \ref{multivariate:thm}.

Write $z_2=(z_i)_{2\le i\le p}$, 
$t_{21}=(t_{i1})_{2\le i\le p}$ and $T_{22}=(t_{ij})_{2\le i,j\le p}$.
Then $z$ and $T$ are represented as block matrices
$$  z = \begin{pmatrix}z_1 \\ z_2 \end{pmatrix}, \qquad
T = \begin{pmatrix} t_{11} & 0 \\ t_{21} & T_{22} \end{pmatrix}. $$
Note that
\begin{eqnarray*}
z'T'T z
&=& (z_1,z_2')
 \begin{pmatrix} t_{11} & t_{21}' \\ 0 & T_{22}' \end{pmatrix}
 \begin{pmatrix} t_{11} & 0 \\ t_{21} & T_{22} \end{pmatrix}
 \begin{pmatrix} z_1 \\ z_2 \end{pmatrix} \\
&=& z_1^2 (t_{11}^2 + t_{21}'t_{21}) + 2 z_1 t_{21}' T_{22} z_2
 + z_2' T_{22}'T_{22} z_2.
\end{eqnarray*}
By taking the expectation with respect to
$t_{11}^2\sim\chi^2_{m+p-1}$, $t_{21}\sim N_{p-1}(0,I_{p-1})$, we have
\begin{eqnarray*}
R_p(z)
&=& z_1^2 (m+p-1+p-1) + R_{p-1}(z_2) \\
&=& z_1^2 (m+2p-2) + R_{p-1}(z_2,\ldots,z_p) \\
&=& \sum_{i=1}^p z_i^2 (m+2p-2i).
\end{eqnarray*} 
Also,
\begin{eqnarray*}
(z'T'T z)^2
&=&  \{ z_1^2 (t_{11}^2 + t_{21}'t_{21}) + 2 z_1 t_{21}' T_{22} z_2
      + z_2' T_{22}'T_{22} z_2 \}^2 \\
&=&  z_1^4 (t_{11}^2 + t_{21}'t_{21})^2 + 4 z_1^2 (t_{21}' T_{22} z_2)^2
    + (z_2' T_{22}'T_{22} z_2)^2 \\
&& +4 z_1^3 (t_{11}^2 + t_{21}'t_{21}) t_{21}' T_{22} z_2 
+ 2 z_1^2 (t_{11}^2 + t_{21}'t_{21}) z_2' T_{22}'T_{22} z_2 \\
&& + 4 z_1 t_{21}' T_{22} z_2 z_2' T_{22}'T_{22} z_2.
\end{eqnarray*}
Noting that
$E[(\chi_\nu^2)^2] = \nu(\nu+2)$, 
we have
\begin{eqnarray*}
S_p(z)
&=& z_1^4 (m+2p-2)(m+2p) + 4 z_1^2 R_{p-1}(z_2) + S_{p-1}(z_2) \\
&& + 2 z_1^2 (m+2p-2) R_{p-1}(z_2) \\
&=& z_1^4 (m+2p-2)(m+2p) + 2 z_1^2 (m+2p) R_{p-1}(z_2) + S_{p-1}(z_2) \\
&=& (z_1^2,\ldots,z_p^2) A_p
 \begin{pmatrix} z_1^2 \\ \vdots \\ z_p^2 \end{pmatrix},
\end{eqnarray*}
where
\allowdisplaybreaks{
\begin{eqnarray*}
A_p
&=&
\left(
\begin{array}{c|c}
 (m+2p)(m+2p-2) & * \\
\hline
\begin{matrix}
 (m+2p)(m+2p-4) \\
 (m+2p)(m+2p-6) \\
 \vdots         \\
 (m+2p) m       
\end{matrix} & A_{p-1}
\end{array}
\right) \\
&=&
\begin{pmatrix}
m+2p   & m+2p   & \cdots & m+2p   \\
m+2p   & m+2p-2 & \cdots & m+2p-2 \\
\vdots & \vdots &        &        \\ 
m+2p   & m+2p-2 &        & m+2
\end{pmatrix}
\odot
\begin{pmatrix}
m+2p-2   & m+2p-4 &        & m      \\
m+2p-4   & m+2p-4 &        & m      \\
         &        &        & \vdots \\ 
m        & m      & \cdots &  m
\end{pmatrix} \\
&=& [ m+2p+2 - 2\min(i,j) ] \odot [ m+2p - 2\max(i,j) ]_{1\le i,j\le p} \\
&=& [ (m+2p+2)(m+2p) - 2(m+2p+2) \max(i,j) - 2(m+2p) \min(i,j) \\
&& + 4 \max(i,j) \min(i,j) ]_{1\le i,j\le p}.
\end{eqnarray*}
Here $\odot$ denotes the elementwise multiplication of matrices.
This means
\begin{eqnarray*}
S_p(z)
&=& (m+2p+2)(m+2p) (\sum_{i=1}^p z_i^2)^2
 - 2(m+2p+2) \sum_{i,j=1}^p \max(i,j) z_i^2 z_j^2 \\
&&
 - 2(m+2p) \sum_{i,j=1}^p \min(i,j) z_i^2 z_j^2
 + 4 (\sum_{i=1}^p i z_i^2)^2.
\end{eqnarray*}
}

\setlength{\bibsep}{1pt}

\end{document}